\documentclass{article}
\usepackage{cite}
\usepackage{amsthm}
\usepackage{graphicx}
\usepackage{amsfonts}

\theoremstyle{plain}
\newtheorem{thm}{Theorem}[section]
\newtheorem{lem}[thm]{Lemma}

\newtheorem{quest}[thm]{Question}

\theoremstyle{definition}
\newtheorem{defn}[thm]{Definition}

\title{Sequences of pseudo-Anosov mapping classes and their asymptotic behavior}
\author{Aaron D Valdivia}

\begin{document}

\maketitle

\section{Introduction}
Consider a surface $S_{g,n}$ with genus $g$ and $n$ punctures with negative Euler characteristic.  The group of self homeomorphisms of the surface up to isotopy is called the \textit{mapping class group}, $\mbox{Mod}^+(S_{g,n})$.  By the Nielsen-Thurston classification \cite{FLP} an element of the mapping class group is either periodic, reducible(i.e. there exists a non-trivial invariant curve), or pseudo-Anosov.  A pseudo-Anosov element is one which fixes a pair of transverse measured singular foliations, 

$$\phi((\mathcal{F}^\pm,\mu^{\pm}))=(\mathcal{F}^{\pm},\lambda^{\pm 1}\mu^{\pm})$$ 

\noindent
up to scaling the measures $\mu^{\pm 1}$ by a constants $\lambda^{\pm 1}$ where $\lambda>1$.  The number $\lambda$ is called the \textit{dilatation} of $\phi$.  The set of dilatations is discrete and therefore bounded away from $1$ \cite{AY} \cite{Iv}.  Furthermore for fixed $g$ and $n$ the minimum, $\delta_{g,n}$, is achieved by some mapping class $\phi\in \mbox{Mod}^+(S_{g,n})$.  One open question about the spectrum of dilataitons is the following one. 

\begin{quest}
What is the value of $\delta_{g,n}$ given $(g,n)$ defining a surface of negative Euler characteristic?
\end{quest}

This question is only answered for a handful of cases of small $g$ and $n$, see \cite{Ko-Los-Song}, \cite{Ham-Song}, \cite{Cho-Ham}.  More is known about the general behavior of these numbers.  Penner explored the behavior for closed surfaces and proved the following theorems.

\begin{thm}\cite{Pen91}\label{lowerbound-thm}
\begin{eqnarray}\label{penner-geq}
\log(\delta_{g,n})\geq\frac{\log(2)}{12g-12+4n}
\end{eqnarray}
\end{thm}

\begin{thm}\cite{Pen91}\label{penner-asymp-thm}
\begin{eqnarray}\label{Penner-eqn}
\log (\delta_{g,0}) \asymp \frac{1}{g}.
\end{eqnarray}
\end{thm}

Tsai \cite{Tsai}(cf \cite{H-K}) continued this investigation for punctured surfaces and showed that for $g=0$ or $g=1$ and $n$ even the behavior is:

$$\log{\delta_{g,n}}\asymp\frac{1}{n}.$$

\noindent
However for surfaces of genus $g>1$ the minimal dilatations behave like, 

$$\log(\delta_{g,n})\asymp\frac{\log(n)}{n}.$$

\noindent
This lead to the following question.

\begin{quest}\cite{Tsai}
What is the behavior for the minimal dilatations for different sequences of $(g,n)$?
\end{quest}

In this paper we provide a partial answer to Tsai's question.

\begin{thm}\label{main-thm}
Given any rational number $r$ the asymptotic behavior along the ray defined by $g=rn$ is, 
$$\log(\delta_{g,n}) \asymp \frac{1}{|\chi(S_{g,n})|},$$
where $\chi(S_{g,n})$ is the Euler characteristic of the surface $S_{g,n}$.
\end{thm}

The proof follows Penner's proof of Theorem \ref{penner-asymp-thm}.  In \cite{Pen91} Penner proves a general lower bound and defines a sequence of pseudo-Anosov mapping classes $\phi_g:S_{g,0}\rightarrow S_{g,0}$ such that $\lambda((\phi_g)^g)$ is bounded by some constant.  We use Penner's lower bound \cite{Pen91} and generalize his examples.  This generalization allows us to construct sequences with dilatation bounded by some constant multiple of $\frac{1}{\mid\chi(S_{g,n})\mid}$.  Given certain choices we can find bounding examples for the upper bound in Theorem \ref{main-thm}.

In Section 2 we recall some known results about pseudo-Anosov mapping classes and train tracks and some techniques of Penner's used in providing the upper bound for closed surfaces.  In Section 3 we define generalized Penner sequences and begin to prove that such a sequence $\phi_n$ has $\log(\lambda(\phi_n))\asymp\frac{1}{\mid\chi(S_{g,n})\mid}$.  In Section 4 we apply our construction and its behavior to the proof of Theorem \ref{main-thm}.

\noindent
\textbf{Acknowledgments:}  I would like to thank Eriko Hironaka for suggesting the problem solved in this paper and for many helpful conversations.  I would also like to thank Beson Farb for reading an earlier version and for his helpful comments.

\section{Background}

In this section we recall some facts about pseudo-Anosov mapping classes and train tracks that will be used later in this paper.  

Dehn showed in \cite{Dehn} that the mapping class group is generated by finitely many Dehn twists $d_x$ where $x$ is a simple closed curve in the surface.  The following theorem of Penner's gives a partial answer to the question of which words of Dehn twists define pseudo-Ansov mapping classes.

\begin{thm}[Penner's Semigroup Criterion]\label{Penner-semigroup-thm} \cite{Pen88}
Suppose $\mathcal{C}$ and $\mathcal{D}$ are each disjointly embedded collections of simple closed curves in an oriented surface $S$.  Suppose $\mathcal{C}$ interesects $\mathcal{D}$ minimally and $\mathcal{C}\cup\mathcal{D}$ fills $S$(i.e. the connected components of $S\slash C\cup D$ have non-negative Euler characteristic).  Let $R(\mathcal{C}^+,\mathcal{D}^-)$ be the semigroup generated by $\{d_c\mid c\in\mathcal{C}\}\cup\{d_d^{-1}\mid d\in\mathcal{D}\}$.  If $\omega\in R(\mathcal{C}^+,\mathcal{D}^-)$, then $\omega$ is pseudo-Anosov if each $d_c$ and $d_d^{-1}$ appear in $\omega$.
\end{thm}

Every pseudo-Anosov mapping class $\phi:S_{g,n}\rightarrow S_{g,n}$ has an associated train track $\tau\subset S_{g,n}$ such that $\phi(\tau)$ is smoothly homotopic into $\tau$, or $\phi(\tau)$ is \textit{carried} by $\tau$.  This homotopy defines a transition matrix $T[a_{i,j}]$ on the edges of $\tau$ where the entry $a_{i,j}$ is the incidence of the edge $i$ with the edge $j$ after applying $\phi$ followed by the homotopy.  This matrix defines a linear action on the vector space of admissable measures of $\tau$.  Each admissable measure of $\tau$ defines a measured singular foliation of $S_{g,n}$ and so finding the real eigenvectors of $T$ is equivalent to finding the invariant foliations where the real eigenvalue is the dilatation.  For a further discussion of train tracks see \cite{P-H} \cite{Pen88}.

Recall that a non-negative matrix $M$ such that $M^n$ is positive for some $n>0$ is said to be Perron-Frobenious.  Such a matrix has a unique real eigenvector whose eigenvalue is the spectral radius of the matrix.  We will use the following lemma to bound the spectral radii of these matrices.  The result is well known and we include the proof for the convenience ofthe reader.

\begin{lem}\label{columnsum-lem}
The spectral radius of a Perron-Frobenius matrix is bounded by the largest column sum.
\end{lem}

\noindent\textbf{Proof:}  Let $M$ be a Perron-Frobenius matrix with spectral radius $\lambda$ and corresponding eigenvector $v$ with norm $1$.  This eigenvector is real and positive.  
$$\mid vM\mid=\lambda=\frac{\sum_{i=1}^nv_im_{ij}}{v_j}\mbox{ for all $j=1...n$}.$$

Let $v_j$ be the largest component of $v$. Then we have:
$$\mid vM\mid=\lambda=\sum_{i=1}^n\frac{v_i}{v_j}m_{ij}.$$
But each term $\frac{v_i}{v_j}\leq 1$ and each term $m_{ij}\geq 0$ and so we have the inequality:
$$\lambda\leq \sum_{i=1}^nm_{ij}.$$
Thus $\lambda$ is bounded by the largest column sum.$\Box$

We will construct pseudo-Anosov maps and corresponding Perron-Frobenius matrices.

As we have already mentioned a pseudo-Anosov mapping class defines a pair of transverse measured singular foliations $(\mathcal{F}^{\pm},\mu^{\pm})$.  The following lemma tells us when a pseudo-Anosov mapping class extends under the forgetful map to another with the same dilatation, see \cite{H-K}. 

\begin{lem}\label{filling-lem}
If $\phi$ is a pseudo-Anosov mapping class on the surface $S_{g,n}$, some subset of the punctures $I$ is fixed setwise, and if none of the points in $I$ are 1-pronged then the punctures may be filled in and the induced mapping class $\tilde{\phi}$ is pseudo-Anosov with $\lambda(\phi)=\lambda(\tilde{\phi})$.
\end{lem} 

An \textit{m-gon} is a possibly punctured disc in $S_{g,n}\backslash\tau$ with $m$ cusps.  We will be able to apply Lemma \ref{filling-lem} with the use of the following Lemma, see \cite{P-H}.

\begin{lem}\label{transitionmatrix-lem} Let $\phi$ be a pseudo-Anosov mapping class and $\tau$ a compatible train track.  Then we have the following:
\begin{enumerate}
\item The singularities and punctures of $S$ defined by the stable foliation of 
$\phi$ are in one-to-one correspondence with the the $m$-gons of $\tau$, and
\item  a singularity or puncture is $m$-pronged if and only if it is contained in an $m$-gon of $\tau$.
\end{enumerate}
\end{lem}

We will use this information in order to construct examples, find transition matrices that bound their dilatation and lastly extend these examples to ones suitable for the $gn$-rays we are interested in.

\section{Penner Sequences}

In this section we will define Penner sequences and give the stepping stones to prove the following theorem.

\begin{thm}\label{sequence-thm}
Given a Penner sequence of mapping classes $\phi_m : F_m \rightarrow F_m$ there is a constant $P$ such that
$$\log(\lambda(\phi_m))\asymp\frac{P}{\mid\chi(F_m)\mid}$$
where $\chi(F_m)$ is the Euler characteristic of the surface $F_m$.
\end{thm}

First we contstruct the surfaces the mapping classes are defined on.  Consider an oriented surface with boundary and punctures, $S_{g,n,b}$, with two sets of disjoint arcs on the boundary components $a^-$ and $a^+$ such that 

$$a^-\cap a^+=\emptyset$$

\noindent
and an orientation reversing homeomorphism, 

$$\iota:a^+\rightarrow a^-.$$  

\noindent
Let $\Sigma_i$ be homeomorphic copies of $S_{g,n,b}$ and let $h_i:S_{g,n,b}\rightarrow \Sigma_i$ be a homeomorphism for each $i\in\mathbb{Z}$.  Set

$$F_{\infty}=\bigcup_{i \in \mathbb{Z}} \Sigma_i/\sim,$$

\noindent
where $y_i \sim y_j$ if, for some $x \in a^+$ and $k\in\mathbb{Z}$,

$$(y_i,y_j)=(h_k(x),h_{k+1}(\iota(x))).$$

The action 

$$\rho=h_{i+1}h^{-1}_i$$ 

\noindent
acts properly discontinuously on $F_{\infty}$.  Then we define

$$F_m=F_\infty\slash\rho^m.$$

We then pick two sets of multicurves $C$ and $D$ on $\Sigma_1$ satisfying Theorem \ref{Penner-semigroup-thm}.  A \textit{connecting curve} is a curve, $\gamma$, on $F_{\infty}$ such that $\gamma\subset\Sigma_1\cup\Sigma_2$, $C\cup\rho(C)\cup\gamma$ is a multicurve, $\gamma$ intersects $D\cup\rho(D)$ minimally, and the set of curves,

$$\displaystyle{J=\{\rho^i(C\cup D\cup\gamma)\}_{-\infty}^{\infty}}$$

\noindent
fills the surface $F_{\infty}$.

\begin{defn}
A sequence of mapping classes $\phi_m:F_m\rightarrow F_m$ is called a Penner sequence if for some $(C,D)$ as in Theorem \ref{Penner-semigroup-thm},
$$\phi_m=\rho d_{\gamma}\omega,$$
where $\omega\in R(C^+,D^-)$ is pseudo-Anosov on $\Sigma_1$.
\end{defn}

Next we want to show that these mapping classes are pseudo-Anosov.  We start with the following lemma about train tracks.  Here we allow our train tracks to have bigons.

\begin{lem}\label{track-lem}
Given a Penner sequence $\phi_m$ there exists invariant train tracks on each surface $F_m$ such that:
\begin{enumerate}
\item  The curves in $C$ and $D$ and the connecting curve are carried on the train track.
\item The images of the curves used to define $\omega$ and the image of the connecting curve under the map $\phi_m$ are carried on the train track.
\end{enumerate}
\end{lem}

\noindent\textbf{Proof:}  We construct a train track on the surface $F_{\infty}$ and then project it to $F_{m}$.  Consider the set of curves $J$ in the definition of $F_{\infty}$.  Then assign positive orientation to all curves,

$$\displaystyle{\{\rho^i(C\cup\gamma)\}_{-\infty}^{\infty}},$$

\noindent
and negative orientation to all curves,

$$\displaystyle{\{\rho^i(D)\}_{-\infty}^{\infty}}.$$

\begin{figure}[h]\label{smooth-fig}
\centering
\includegraphics[scale=.5]{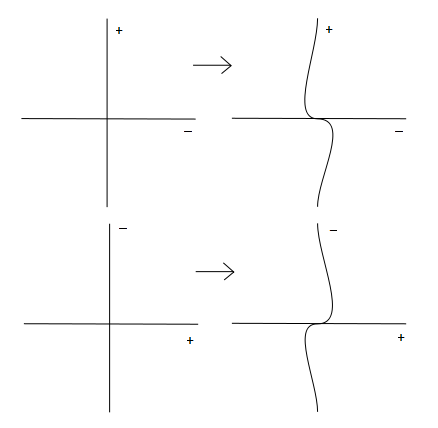}
\caption{Smoothing}
\end{figure}

\noindent
We then smooth the intersections according to Figure \ref{smooth-fig}.

This provides a train track, $\tau$, which may be projected to a track $\tau_m$ on $F_m$.  It is easy to see that all the curves in $J$ are carried by the train track.  If we consider a curve $x$ that is carried by $\tau_m$ and perform a Dehn twist by an element $y\in J$ then the resulting curve $d_y(x)$ is carried by the edges carrying $x$ and the edges that carry the element $y$.  Since the train track is symmetric with respect to the rotation map we are done.$\Box$

The next lemma allows us to compute the entries of a transition matrix on this train track.  

\begin{lem}\label{measure-lem}
Consider $x,y\in J$ where, 

$$J=\displaystyle{\cup_{i=1...m}\rho^i(C\cup D\cup \gamma)}.$$  

If $\mu_x$ and $\mu_y$ are the elements of $W_{\tau_m}$ induced by the curves $x$ and $y$ and $d_x^{\star}$ is the map on $W_{\tau_m}$ induced by $d_x$ then 
$$d_x^{\star}(\mu_y)=\mu_y+i(x,y)\mu_x$$
\end{lem}

\noindent\textbf{Proof:}  Collapsing a curve is invariant up to homotopy so the edges that elements of $J$ collapse to on the train track are well defined.  Therefore we can consider a curve as a free homotopy class.  If two curves $x,y\in J$ are in minimal position then performing a Dehn twist on $y$ about $x$ will send $y$ to $d_x(y)=y+i(x,y)x$.  We collapse the homology elements to $\tau_m$ and obtain the desired result.$\Box$

\begin{lem}\label{PF-lem}
The matrix $T_{\phi_m}$ is Perron-Frobenius for each $m \geq 2$.
\end{lem}

\noindent\textbf{Proof:}  This follows from the previous lemma and the fact that given the mapping class $\phi_m^m$ where we perform Dehn twists about all the curves in, 

$$J_m=\displaystyle{\cup_{i=1...m}\rho^i(C\cup D\cup \gamma)},$$ 

\noindent
each curve is connected to another through at most $m(r+s+1)$ curves.  Therefore $T_{\phi_m}^{m^2(r+s+1)}$ is strictly positive where $r=\sharp C$ and $s=\sharp D$.$\Box$

\bigskip

\noindent\textbf{Remark:}  An admissable measure on a bigon track defines a measured foliation up to an equivalence of the admissable measures \cite{P-H}.  We only consider a subset of the admissable measures when computing transition matrices.  Since we have show that the transition matrix on these measures is Perron-Frobenius the eigenvector is positive, defining an invariant foliation for the mapping class.  We will prove that these mapping classes are pseudo-Anosov in the next section, therefore the invariant expanding foliation is unique and we need not worry about the equivalent measures or measures outside the considered subset.

\section{Asymptotic Behavior}

\noindent\textbf{Proof of Theorem \ref{sequence-thm}:}  The mapping classes $(\phi_m)^m$ are pseudo-Anosov by Penner's semigroup criteria and so the mapping classes $\phi_m$ are as well.  Lemma \ref{track-lem} gives an invariant bigon train track.

As stated in Theorem \ref{columnsum-lem} the spectral radius of a Perron-Frobenius matrix is bounded by the largest column sum of the matrix.  So now we would like to compute the matrix defining the action on the transverse measures.  This matrix will be Perron-Frobenius by Lemma \ref{PF-lem} and can be computed using Lemma \ref{measure-lem}.

The map $\rho$ permutes the curves of $J_m=\displaystyle{\cup_{i=1...m}\rho^i(C\cup D\cup \gamma)}$ and so the induced map on the space of weights spaned by $\mu_1...\mu_{m(r+s+1)}$, where $r=\sharp C$ and $s=\sharp D$, is defined by a block permutation matrix.

$$M_{\rho}=\left(
\begin{array}{ccccc}
0&I&0&.&0\\
0&0&I&.&0\\
.&.&.&.&.\\
I&0&0&.&0
\end{array}\right)$$

\noindent
The Dehn twist about the connecting curve gives the map with transition matrix defined by,

$$M_{d_{\gamma}}=\left(
\begin{array}{ccccc}
U&0&0&.&0\\
V&I&0&.&0\\
.&.&.&.&.\\
0&0&0&.&I
\end{array}\right).$$

\noindent
Last the transition matrix for the map induced by the word $\omega$ is given below.

$$M_{\omega}=\left(
\begin{array}{ccccc}
W&0&0&.&0\\
0&I&0&.&0\\
0&0&I&.&0\\
.&.&.&.&.\\
X&0&0&.&I
\end{array}\right)$$

\noindent
Then the matrix for the map $\phi_m$ is given below by matrix multiplication after making the identifications $WU=Y$ and $XK=Z$.

$$M_{\phi_m}=\left(\begin{array}{ccccccc}
0&Y&0&0&.&.&0\\
0&V&I&0&.&.&0\\
0&0&0&I&.&.&0\\
.&.&.&.&.&.&.\\
.&.&.&.&.&.&.\\
0&0&0&0&.&.&I\\
I&Z&0&0&.&.&0
\end{array}\right)$$

The matrix $M$ is an $m\times m$ block matrix of $r+s+1\times r+s+1$ blocks.  The matrix $Y$ depends on the word $\omega$.  The matrix $V$ may have non-zero entries in the last column except that the last row must be zero since $\gamma\cap\rho(\gamma)=\emptyset$.  The observation $V^2=0$ will be important later.

$$V=\left(\begin{array}{cccccc}
0&.&.&.&0&v_1\\
0&.&.&.&0&v_2\\
.&.&.&.&.&.\\
.&.&.&.&.&.\\
0&.&.&.&0&v_{rs}\\
0&.&.&.&0&0
\end{array}\right)$$

\noindent
The matrix $Z$ will depend on $\omega$ as well but only has non-zero entries in the last row.  Now we want to consider the matrix $M^m$.  Inductively we see that for $1<k<m$ the matrix $M^k$ is given by the following matrix.  Here we use the fact that $V^2$ is the zero matrix.

$$M_{\phi_m}^k=\left(
\begin{array}{cccccccccc}
0 &0&    0&    .& YV&   Y& 0& 0& .&0\\
.& .&    .&    .& 0&    V& I& 0& .&.\\
.& .&    .&    .& .&    0& 0& I& .&0\\
0& 0&    .&    .& .&    .& .& 0& .&I\\
I& Z&    0&    .& .&    .& .& .& .&0\\
0& Y+ZV& Z&    .& .&    .& .& .& .&.\\
.& YV&   Y+ZV& .& .&    .& .& .& .&.\\
.& 0&    YV&   .& .&    .& .& .& .&.\\
.& .&    0&    .& 0&    .& .& .& .&.\\
.& .&    .&    .& C&    0& .& .& .&.\\
0& 0&    0&    .& Y+ZV& Z& 0& 0& .&0
\end{array}\right)$$

\noindent
Then we can find the transition matrix for the $m$th iterate.

$$M^m=\left(
\begin{array}{ccccccc}
Y& YZ  &   0   &    .& .& 0   &    YV\\
V& Y+ZV&   Z   &    .& .& .   &    0\\
0& YV  &   Y+ZV&    .& .& .   &    .\\
.& 0   &   YV  &    .& .& .   &    .\\
.& .   &   0   &    .& .& 0   &    .\\
.& .   &   .   &    .& .& Z   &    0\\
0& .   &   .   &    .& .& Y+ZV&    Z\\
C& 0   &   0   &    .& .& YV  &   Y+ZV
\end{array}\right)$$

Then this matrix is Perron-Frobenius by Lemma \ref{PF-lem} and by Theorem \ref{columnsum-lem} the spectral radius is bounded by the largest column sum.  A block column sum is either equal to a column sum of $YV+Z+ZV$, $YZ+Y+ZV+YV$, or $Y+V+Z$.  Therefore the dilatation of the $m$th iterate is bounded by a constant, say $P$.  This tells us that 

$$\log(\lambda(\phi_m))\leq \frac{P}{m}.$$

\noindent
Then Theorem \ref{lowerbound-thm} with the upper bound just given finishes the proof.$\Box$

\medskip

\noindent
Next we use this to prove Theorem \ref{main-thm}

\medskip

\begin{figure}\label{proof-fig1}
\centering
\includegraphics[scale=.5]{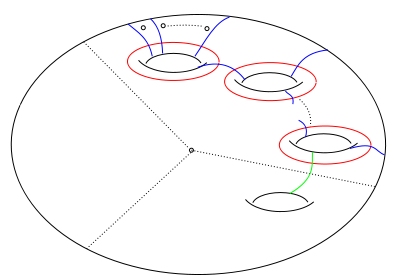}
\caption{Choice for curves}
\end{figure}

\begin{figure}\label{proof-fig2}
\centering
\includegraphics[scale=.35]{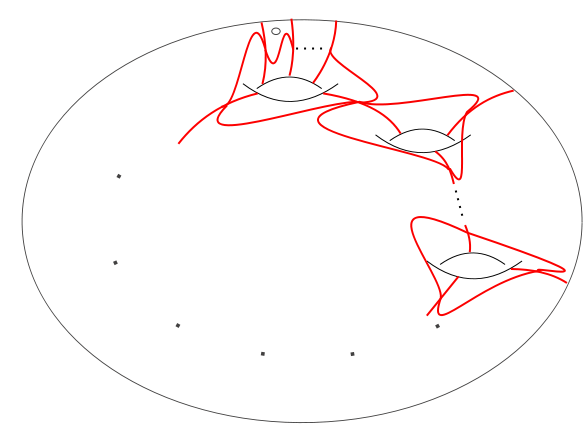}
\caption{Train track}
\end{figure}

\noindent\textbf{Proof of Theorem \ref{main-thm}:}  With Penner's lower bound we only need the upper bound to prove the asymptotic behavior.  Suppose a ray has slope $\frac{p}{q}$ with $(p,q)=1$ then let $S_{g,n,1}$ have $g=p$, and $n=q$ and create any Penner sequence with $a^+$ being an arc on the boundary component and $a^-$ another disjoint arc on the same boundary component.  This sequence of mapping classes then has the required upper bound on dilatation.  This gives sequences on $gn$-rays through $(2,0)$.  Further if we choose our curves as in Figure \ref{proof-fig1}, which is shown with a chosen connecting curve as well, then we can find a train track for $\phi_m$ given in Figure \ref{proof-fig2}.  

From this we can see by Theorem \ref{transitionmatrix-lem} that the two fixed punctures are not 1-pronged.  Filling in both fixed punctures we obtain sequences of mapping classes with the same dilatation and two fewer punctures, the sequences for the lines passing through the origin.$\Box$

\bibliography{mybib}
\bibliographystyle{alpha}
\end{document}